\documentclass[12pt,reqno,twoside]{amsart}
\usepackage{amsmath}
\usepackage{amssymb}
\usepackage{amsthm}
\newtheorem {theorem}    {Theorem}[section]
\newtheorem {lemma}      [theorem]    {Lemma}

\newtheorem {corollary}  [theorem]    {Corollary}
\newtheorem {proposition}[theorem]    {Proposition}

\theoremstyle{definition}

\newcounter{AbcT}

\numberwithin{equation}{section}

\newcommand {\R} {{\mathbb R}}

\newcommand {\Z} {{\mathbb Z}}

\renewcommand{\liminf}{\varliminf}

\DeclareMathOperator{\supp}{supp}

 \DeclareMathOperator{\SL}{SL}

\DeclareMathOperator{\GL}{GL}

\newcommand {\IGNORE}[1]  {}

\newcommand{\nope} {}

\ifx\pdfoutput\undefined \usepackage[hypertex]{hyperref} \else \usepackage[pdftex,pdfstartview=FitH]{hyperref} \fi

\begin{document}

\title[$\SL (2, \R)$-invariant measures] { Ratner's theorem on $\SL (2, \R)$-invariant measures}

\begin{abstract}We give a relatively short and self contained proof
of Ratner's theorem in the special case of $\SL (2, \R)$-invariant
measures. \end{abstract}

\author {Manfred Einsiedler}
\address{Mathematics Department, The Ohio State University, 231 W. 18th Avenue, Columbus, Ohio 43210}
\thanks{The author acknowledges support of NSF Grant 0509350.
This research was partially conducted while the author
was employed by the Clay Mathematics Institute as a Research
Scholar.}

\date{\today}

\maketitle

\section {Introduction}

M.~Ratner proved in a series of papers
\cite{Ratner-invent-solvable,Ratner-acta-measure,Ratner-measure-rigidity,Ratner91t,Ratner-bull-distribution}
very strong results on invariant measures and orbit closures for
certain subgroups $H$ of a Lie group $G$ --- where $H$ acts on the
right of $X = \Gamma \backslash G$ and $\Gamma< G$ is a lattice.
More concretely, $H$ needs to be generated by one parameter
unipotent subgroups, and the statements all are of the form that
invariant measures and orbit closures are always algebraic as
conjectured earlier by Raghunathan. Today these theorems are applied
in many different areas of mathematics.

While there are some very special cases for which the proof
simplifies \cite{Ratner-sl2,Witte}, the general proof requires a
deep understanding of the structure of Lie groups and ergodic
theory. The aim of this paper is to give a self-contained more
accessible proof of the classification of invariant and ergodic
measures for subgroups $H$ isomorphic to $\SL (2, \R)$.  While this
is as well a special case of Ratner's theorem, it is a rich class
since $G$ can be much larger than $H$. Moreover, the proof for this
class is more accessible in terms of its requirements.  The methods
used are not new, in particular most appear also in earlier work
\cite{Ratner-joinings, Ratner-acta-measure} of Ratner, but it does
not seem to be known that the following theorem allows also a
relatively simple and short proof.

\begin{theorem}  \label{Main theorem}
Let $G$ be a Lie group, $\Gamma< G$ a discrete subgroup, and $H< G$
a subgroup isomorphic to $\SL (2, \R)$.  Then any $H$-invariant and
ergodic probability measure $\mu$ on $X = \Gamma \backslash G$ is
homogeneous, i.e.\ there exists a closed connected subgroup $L < G$
containing $H$ such that $\mu$ is $L$-invariant and some $x _ 0 \in
X$ such that the $L$-orbit $x_0 L$ is closed and supports $\mu$.  In
other words $\mu$ is an $L$-invariant volume measure on $x_0L$.
\end{theorem}

As we will see a  graduate student, who started to learn or is
willing to learn the very basics of Lie groups and ergodic theory,
should be able to follow the argument (no knowledge of radicals or
structure theory of Lie groups and no knowledge of entropy is
necessary). We will discuss the requirements in
Section~\ref{ingredients}. The initiated reader will notice that
this approach generalizes without too much work to other semisimple
groups $H<G$ without compact factors.  However, to keep the idea
simple and to avoid unnecessary technicalities we only treat the
above case.

In the next section we give some motivation for the above and
related questions.  In particular, we will discuss an application of
Ratner's theorems where the above special case is sufficient. This
paper is best described as an introduction to Ratner's theorem on
invariant measures, and is not a comprehensive survey of this area
of research. The reader seeking such a survey is referred to
\cite{survey-Kleinbock} and for some more recent developments to
\cite{Einsiedler-Lindenstrauss-survey}.

The author would like to thank M.~Ratner for comments on an earlier
draft of this paper.

\section {Motivation}

Let $G$ be a closed linear group, i.e.\ a closed subgroup of $\GL(n,\R)$ (with that assumption some definitions are easier to make).
Let $\Gamma<G$ be a discrete subgroup so that $X=\Gamma\backslash G$ is a locally compact space with a natural $G$-action:
\[
\mbox{for }g\in G,x\in X\mbox{ let }g.x=xg^{-1}.
\]
Now let $H<G$ be a closed subgroup and restrict the above action to the subgroup $H$. The question about the properties of the resulting $H$-action has many connections to various a priori non-dynamical mathematical problems, and is from that point of view but also in its own light highly interesting.

The most basic question, vaguely formulated, is how $H$-orbits $H.x_0$ for various points $x_0\in X$ look like. Here one can make restrictions on $x_0$ or not, and ask, more precisely, either about the distribution properties of the orbit or about the nature of the closure of the orbit.

If $X$ carries a $G$-invariant probability measure $m_X$, which in
many situations is the case, then $m_X$ is called the {\em Haar
measure} of $X$ and $\Gamma$ is by definition a {\em lattice}. In
that situation one can restrict these questions to $m_X$-typical
points and by doing so one has entered the realm of ergodic theory.
Rephrased the basic question is now whether $m_X$ is $H$-ergodic.
Recall that by definition $m_X$ is $H$-ergodic if every measurable
$f$ that is $H$-invariant is constant a.e.\ with respect to $m_X$,
and note that while $H$-invariance of $m_X$ is inherited from
$G$-invariance, the same is not true for $H$-ergodicity. The
characterization of $H$-ergodicity in this context has been given in
varying degrees of generality by many authors mostly before 1980,
see \cite[Chpt.\ 2]{survey-Kleinbock} for a detailed account.  The
power of this characterization is that often --- unless there are
obvious reasons for failure of ergodicity --- the Haar measure turns
out to be $H$-ergodic. Moreover, assuming ergodicity one of the most
basic theorems in abstract ergodic theory, the ergodic theorem,
states that a.e.\ points equidistribute in $X$. (For the notion of
equidistribution we need that $H$ is an amenable subgroup. The fact
that $ \SL (2, \R)$ is not amenable actually makes it harder to
apply Theorem \ref{Main theorem} --- and is the reason why the proof
simplifies.) In particular, a.e.\ $H$-orbit is dense. This can be
seen as the first answer to our original question.

Let us assume from now on that $\Gamma$ is a lattice (some of what follows holds more generally for any discrete $\Gamma$ but not all). The above discussion around $H$-ergodicity is only the first step in understanding the structure of $H$-orbits. In general, there is no reason to believe that a similarly simple answer is possible for all points of $X$.  This is especially true for general dynamical systems but as we will discuss also in the algebraic setting we consider here.  More surprisingly, there are cases where we can understand {\bf all} $H$-orbits respectively believe that it is possible to understand {\bf all} $H$-orbits.
To be able to describe this we need to recall a few notions: $u \in G$ is {\em unipotent} if $1$ is the only eigenvalue of the matrix $u$, $a \in G$ is $\R$-{\em diagonalizable} if it is diagonalizable as a matrix over $\R$.  A subgroup $U < G$ is a {\em one-parameter unipotent subgroup} if $U$ is the image of a homomorphism $t \in \R \mapsto u _ t \in U$ with $u _ t$ unipotent for all $t \in \R$.  E.g.\ in $G = \SL (n, \R)$ the subgroup
\begin{equation*}
U = \begin{pmatrix}1 & * &\cdots & *\\ 0&1&\ddots&\vdots \\ \vdots&\ddots&\ddots&*
\\ 0 &\cdots & 0 &1
\end{pmatrix}
\end{equation*}
contains only unipotent elements and is generated by one-parameter unipotent subgroups.  A subgroup $A < G$ is $\R$-{\em diagonalizable} if for some $g \in \GL (n, \R)$ the conjugate $g A g ^{- 1}$ is a subgroup of the diagonal subgroup.

\subsection{Unipotently generated subgroups and Oppenheim's conjecture}
One case of subgroups we today understand quite well is when $H$ is
generated by one-parameter unipotent subgroups (as it is the case
for $\SL (2, \R)$). As mentioned before this is thanks to the
theorem of M.~Ratner \cite{Ratner91t} which says that all $H$-orbits
are well-behaved as conjectured by Raghunathan earlier. Every
$H$-orbit $H.x_0$ is dense in the closed orbit $L.x_0$ of a closed
connected group $L>H$, and the latter orbit $L.x_0$ supports a
finite $L$-invariant volume measure $m_{L.x_0}$. If $H$ is itself a
unipotent one-parameter subgroup the orbit $H.x_0$ is
equidistributed with respect to this measure $m_{L.x_0}$. These
theorems were later extended by Ratner \cite{Ratner-padic} herself,
as well as by Margulis and Tomanov \cite{Margulis-Tomanov} to the
more general setting of products of linear algebraic groups over
various local fields ($S$-algebraic groups).

A bit more technical is the following question: What are the
$H$-invariant probability measures? Here it suffices to restrict to
$H$-invariant and ergodic measures --- the general theorem of the
ergodic decomposition states that any $H$-invariant probability
measure can be obtained by averaging $H$-invariant and ergodic
measures. Therefore, if we understand the latter measures we
understand all of them. Ratner showed that all $H$-invariant and
ergodic probability measures are of the form $m_{L.x_0}$ as
discussed above --- Theorem \ref{Main theorem} is a special case of
this. As it turns out this question is crucial for the proof of the
topological theorem regarding orbit closures mentioned above. Namely
using her theorem on measure classification Ratner then proves first
that the orbit of a unipotent group always equidistributes with
respect to an ergodic measure, and finally uses this to prove her
topological theorem
--- this is similar to the earlier discussion of $m_X$-typical
points. However, these steps are quite involved: First of all it is
not clear that a limit distribution coming from the orbit of a
one-parameter unipotent subgroup is a probability measure since the
space might not be compact. However, earlier work of Dani
\cite{Dani-inventiones, Dani-I, Dani-II} (which extend work by
Margulis \cite{Margulis-nondivergence}) shows precisely this. Then
it is not clear why such a limit is ergodic and why it is
independent of the times used in the converging subsequence ---
without going into details let us just say that the proof relies
heavily on the structure of the ergodic measures and the properties
of unipotent subgroup.

Before Ratner classified in her work all orbit closures Margulis
used a special case of this to prove Oppenheim's conjecture, which
by that time was a long standing open conjecture. This conjecture
concerns the values of an indeterminate irrational quadratic form
$Q$ in $n$ variables at the integer lattice $\Z ^n$, and Margulis
theorem \cite{Margulis-Oppenheim-proof} states that under these
assumptions $Q(\Z^n)$ is dense in $\R$ if $n\geq 3$. (It is not hard
to see that all of these assumptions including $n\geq 3$ are
necessary for the density conclusion.) The proof consists of
analyzing all possible orbits of $\operatorname{SO}(2,1)$ on
$X_3=\SL(3,\Z) \backslash \SL(3,\R)$. Even though
$\operatorname{SO}(2,1)$ is essentially a quotient of $\SL (2, \R)$
Theorem \ref{Main theorem} does not imply immediately Oppenheim's
conjecture since for the non-amenable group $\SL (2, \R)$ it is not
obvious how to find an invariant measure on the closure of an orbit.

\subsection{Diagonalizable subgroups}

The opposite extreme to unipotent elements are $\R$-diagonalizable elements, so it is natural to ask next about the case of a $\R$-diagonalizable subgroups $A<G$: What do the closures of $A$-orbits look like? What are the $A$-invariant and ergodic probability measures? As we will see this case is more difficult in various ways, in particular we do not have complete answers to these questions.

{\bf Rank one:}
Let $G = \SL (2, \R)$, $\Gamma = \SL (2, \Z)$, and set $X _ 2 = \Gamma \backslash G$.  Let $A< \SL (2, \R) $
be the diagonal subgroup.  Then the action of $A$ on $X _ 2$ can also be described as the geodesic flow on the unit tangent bundle of the modular surface.  Up to the fact that the underlying space $X _ 2 $ is not compact this is a very good example of an Anosov flow.  The corresponding theory can be used to show that there is a huge variety of orbit closures and $A$-invariant ergodic probability measures. So the answer is in that case that there is no simple answer to our questions --- but at least we know that.

{\bf Higher rank:} We replace `2' by `3' and encounter very
different behaviour. Let $G = \SL (3, \R)$, $\Gamma = \SL (3, \Z)$,
and set $X _ 3 = \Gamma \backslash G$.  Let $A< \SL (2, \R) $ be the
diagonal subgroup, which this time up to finite index is isomorphic
to $\R ^2$. Margulis, Furstenberg, Katok, and Spatzier conjectured
that for the action of $A$ on $X_3$ there are very few $A$-invariant
and ergodic probability measures, in particular that they again are
all of the type $m_{L.x_0}$ for some $L>A$ and $x_0 \in X_3$ with
closed orbit $L.x_0$. One motivation for that conjecture is
Furstenberg's theorem \cite{Furstenberg-67} on $ \times 2$, $\times
3$-invariant closed subsets of $\R / \Z$ which states that all such
sets are either finite unions of rational points or the whole space.
This can be seen as an abelian analogue to the above problem.

For orbit closures the situation is a bit more complicated: $G$ contains an isomorphic copy $L$ of the subgroup $\GL(2,\R)$ embedded into the upper left 2-by-2 block (where the lower right entry is used to fix the determinant). Now $\Gamma L$ is a closed orbit of $L$ and the $A$-action inside this orbit consists of the rank one action discussed above and one extra direction which moves everything towards infinity. (This behaviour of the $L$-orbit and the $A$-orbits needs of course justification, which in this case can be done algebraically.) Therefore, in terms of orbit closures the situation for $A$-orbits inside this $L$-orbit is as bad as for the corresponding action on $X_2$. However, it is possible to avoid this issue and to formulate a meaningful topological conjecture: Margulis conjectured that all bounded $A$-orbits on $X_3$ are compact, i.e.\ the only bounded orbits are periodic orbits for the $A$-action (which in fact all arise from a number theoretic construction).

Margulis \cite{Margulis-Oppenheim-conjecture} also noted that the
question regarding orbits in this setting is related to a long
standing conjecture by Littlewood. Littlewood conjectured around
1930 that for any two real numbers $\alpha, \beta \in \R$ the vector
$(\alpha, \beta)$ is well approximable by rational vectors in the
following multiplicative manner:
\begin{equation*}
\liminf _{n \rightarrow \infty} n \|n \alpha\|\|n \beta\|=0,
\end{equation*}
where $\|u\|$ denotes the distance of a real number $u \in \R$ to $\Z$. Here $n$ is the common denominator of the components of the rational vector that approximates $(\alpha, \beta)$, and instead of taking the maximum of the differences along the $x$-axis and the $y$-axis we instead measure the quality of approximation by taking the product of the differences. The corresponding dynamical conjecture states that certain points (defined in terms of the vectors $(\alpha, \beta)$) all have unbounded orbit (where actually only a quarter of the acting group $A$ is used).

Building on earlier work of E.~Lindenstrauss
\cite{Lindenstrauss-Quantum} and a joint work of the author with
A.~Katok \cite{Einsiedler-Katok-02} we have obtained together
\cite{Einsiedler-Katok-Lindenstrauss} a partial answer to the
conjecture on $A$-invariant and ergodic probability measures: If the
measure has in addition positive entropy for some element of the
action, then it must be the Haar measure $m_{X_3}$ --- this
generalizes earlier work of Katok, Spatzier, and Kalinin
\cite{Katok-Spatzier-96, Katok-Spatzier-98, Kalinin-Spatzier-02},
and related work by Lyons \cite{Lyons-88}, Rudolph
\cite{Rudolph-90}, and Johnson \cite{Johnson-2-3} in the abelian
setting of $\times 2, \times 3$. For Littlewood's conjecture we show
also in \cite{Einsiedler-Katok-Lindenstrauss} that the exceptions
form at most a set of Hausdorff dimension zero. Roughly speaking the
classification of all $A$-invariant probability measures with
positive entropy can be used to show that very few $A$-orbits can
stay within a compact subset of $X_3$, which by the mentioned
dynamical formulation of Littlewood's conjecture is what is needed.

Most of the proof of this theorem consists in showing that positive
entropy implies invariance of the measure $\mu$ under some subgroup
$H<\SL (3, \R)$ that is generated by one-parameter unipotent
subgroups.  Then one can apply Ratner's classification of invariant
measures to the $H$-ergodic components of  $\mu$. However, in this
case (unless we are in the easy case of $H=\SL (3, \R)$) the
subgroup $H$ is actually isomorphic to $\SL (2, \R)$.  Therefore,
Theorem \ref{Main theorem} is sufficient to analyze the $H$-ergodic
components.

The more general case of the maximal diagonal subgroup $A$ acting on
$X_n$ for $n\geq 3$ is also treated in
\cite{Einsiedler-Katok-Lindenstrauss} (always assuming positive
entropy). Even more generally one can ask about any
$\R$-diagonalizable subgroup of an (algebraic) linear group.
However, here there are unsolved technical difficulties that prevent
so far a complete generalization. Joint ongoing work
\cite{Einsiedler-Lindenstrauss-maximally-split} of the author with
E.~Lindenstrauss solves these problems for maximally
$\R$-diagonalizable subgroups (more technically speaking, for
maximal $\R$-split tori $A$ in algebraic groups $G$ over $\R$ and
similarly also for $S$-algebraic groups). This approach uses results
from \cite{Einsiedler-Katok-nonsplit} and
\cite{Einsiedler-Lindenstrauss-low-entropy}. For a more complete
overview of these results and related applications see the survey
\cite{Einsiedler-Lindenstrauss-survey}.

\section{Ingredients of the proof} \label{ingredients}

We list the facts and notions needed for the proof of Theorem \ref{Main theorem}, all of which, except for the last one, can be found in any introduction to Lie groups resp.\ ergodic theory.

\subsection{ The Lie group and its Lie algebra}
The Lie algebra $\mathfrak g$ of $G$ is the tangent space to $G$ at the identity element $e \in G$.
The exponential map $\exp: \mathfrak g \rightarrow G$ and the locally defined inverse, the logarithm map, give local isomorphisms between $\mathfrak g$ and $G$.  For any $g \in G$ the derivative of the
conjugation map is the adjoint transformation $\operatorname{Ad} _ g: \mathfrak g \rightarrow \mathfrak g$ and satisfies $\exp \operatorname{Ad} _ g (v) = g \exp (v) g ^ {- 1}$ for $g \in G$ and $v \in \mathfrak g$. For linear groups this could not be easier, the Lie algebra is a linear subspace of the space of matrices, $\exp(\cdot)$ and $\log(\cdot)$ are defined as usual by power series, and the adjoint transformation $\operatorname{Ad} _ g$ is still conjugation by $g$.

Closed subgroups $L < G$ are almost completely described by their respective Lie algebras $\mathfrak l$ inside $\mathfrak g$ as follows. Let $L ^\circ$ be the connected component of $L$ (that contains the identity $e$).  Then the Lie algebra $\mathfrak l$ of $L$ (and $L ^\circ$) uniquely determines $L ^\circ$ --- $L ^\circ$ is the subgroup generated by $\exp (\mathfrak l)$.  (Moreover, any element $\ell \in L$ sufficiently close to $e$ is actually in $L^\circ$ and equals $\ell= \exp (v)$ for some small $v \in \mathfrak l$.)

Using an inner product on $\mathfrak g$ we can define a left invariant Riemannian metric $d(\cdot,\cdot)$ on $G$.  We will be using the restriction of $d(\cdot,\cdot)$ to subgroups $L < G$ and denote by $B_r^L$ the $r$-ball in $L$ around $e \in L$.

If $\Gamma <G$ is a discrete subgroup, then $X=\Gamma \backslash G$ has a natural topology and in fact a metric defined by $d(\Gamma g, \Gamma h)=\min_ {\gamma \in \Gamma}d(g, \gamma h)$ for any $g,h \in G$ (which uses left invariance of $d(\cdot,\cdot)$). With this metric and topology $X$ can locally be described by $G$ as follows. For any $x \in X$ there is an  $r>0$ such that the map $\imath: g\mapsto xg$ is an homeomorphism between $B_r^G$ and a neighborhood of $y$. Moreover, if $r$ is small enough $\imath:B_r^G \rightarrow X$ is in fact an isometric embedding. For a given $x$ a number $r>0$ with these properties we call an {\em injectivity radius at $x$}.

\subsection{Complete reducibility and the irreducible representations of $ \SL (2, \R)$} \label{Representations}
The first property of $\SL (2, \R)$ we will need is the following standard fact.  Let $V$ be a finite dimensional real vector space and suppose $\SL (2, \R)$ acts on $V$.  Then any $\SL (2, \R)$-invariant subspace $W < V$ has an $\SL (2, \R)$-invariant complement $W ' < V$ with $V = W \oplus W '$.

The above implies that all finite dimensional representations of $\SL (2, \R)$ can be written as a direct sum of irreducible representations.  The second fact we need is the description of these irreducible representations. Let $A= \begin{pmatrix} 1 \\ 0 \end{pmatrix}$ and $B = \begin{pmatrix}0 \\ 1 \end{pmatrix}$ denotes the standard basis of $\R ^ 2 $ so that $\begin{pmatrix}1 & t \\ 0 & 1 \end{pmatrix} A = A$ and  $\begin{pmatrix}1 & t \\ 0 & 1 \end{pmatrix} B = B+t A$.  Any irreducible representation  is obtained as a symmetric tensor product $\operatorname{Sym}^n(\R ^ 2)$ of the standard representation on $\R ^ 2$ for some $n$.  $\operatorname{Sym}^n(\R ^ 2)$ has $A^n,A^{n-1}B,\ldots,B^n$ as a basis, and every element we can view as a homogeneous polynomial $p(A,B)$ of degree $n$.  The action of  $\begin{pmatrix}1 & t \\ 0 & 1 \end{pmatrix}$ can now be described by substitution, $p(A,B)$ is mapped to $p(A,B+tA)$.  More concretely, $p(A,B)=c_0A^n+c_{1}A^{n-1}B+\cdots+c_nB^n$  is mapped to
\begin{align*}
p(A,B+tA)=&(c_0+c_{1}t+\cdots+c_nt^n)
A^n+\\&(c_{1}+\cdots+c_nnt^{n-1})A^{n-1}B+\\&\cdots+c_nB^n,
\end{align*}
where the coefficients in front of the various powers of $t$ are the original components of the vector $p(A,B)$ multiplied by binomial coefficients.  Notice that all components of $p(A,B)$ appear in the image vector in the component corresponding to $A^n$. Moreover, for any component of $p(A,B)$ the highest power of $t$ it gets multiplied by appears in the resulting component corresponding to $A^n$. For that reason, when $t$ grows (and say $p (A,B)$ is not just a multiple of $A ^ n$) the image of $p(A,B)$ under $\begin{pmatrix}1 & t \\ 0 & 1 \end{pmatrix}$ will always grow fastest in the direction of $A ^ n$ when $t\rightarrow\infty$.

\subsection{Recurrence and the ergodic theorem}
Let $ (X, \mu)$ be a probability space, and let $T:X \rightarrow X$ be measure preserving. Then for any set $B \subset X$ of positive measure and a.e.\ $x \in B$ there are infinitely many $n$ with $T^nx \in B$ by Poincar\'e recurrence.

Now suppose $u _ t: X \rightarrow X$ for $t \in \R$ is a one parameter flow acting on $X$ such that $\mu$ is $u _ \R$-invariant and ergodic.  Then for any $f \in L^1 (X, \mu)$ and $\mu$-a.e.\ $x \in X$ we have
\begin{equation*}
\frac{1}{T} \int_ 0 ^T f(u_t(x))\operatorname{d}\! t \rightarrow \int_ X f \operatorname{d}\mu\mbox{ for }T \rightarrow \infty
\end{equation*}
This is Birkhoff's pointwise ergodic theorem for flows.

\subsection{Mautner's phenomenon for $\SL(2,\R)$}
To be able to apply the ergodic theorem as stated in the last
section in the proof of Theorem~\ref{Main theorem} we will need to
know that the $\SL(2, \R)$-invariant and ergodic probability measure
is also ergodic under a one-parameter flow. The corresponding fact
is best formulated in terms of unitary representations and is due to
Moore \cite{Moore} and is known as the {\em Mautner phenomenon}. For
completeness we prove the special case needed.

\begin{proposition} \label{Mautner}
Let $\mathfrak H$ be a Hilbert space, and suppose $\phi:\SL(2,\R)
\rightarrow U(\mathfrak H)$ is a continuous unitary representation
on $\mathfrak H$. In other words, $\phi$ is a homomorphism into the
group of unitary automorphisms $\mathbb{U}(\mathfrak H)$ of
$\mathfrak H$ such that for every $v \in \mathfrak H$ the vector
$\phi(g)(v) \in \mathfrak H$ depends continuously on $g \in
\SL(2,\R)$. Then any vector $v \in \mathfrak H$ that is invariant
under the upper unipotent matrix group $U=\left\{\begin{pmatrix} 1 &
* \\ 0 & 1\end{pmatrix}\right\}$ is in fact invariant under $\SL
(2,\R)$.
\end{proposition}

Since any measure preserving action on $(X,\mu)$ gives rise to a
continuous unitary representation on $\mathfrak H=L^2(X,\mu)$ the
above gives immediately what we need (see also
\cite[Prop.~5.2]{Ratner-acta-measure} for another elementary
treatment):

\begin{corollary} \label{Mautner ergodicity}
Let $\mu$ be an $H$-invariant and ergodic probability measure on
$X=\Gamma\backslash G$ with $\Gamma<G$ discrete, and $H<G$
isomorphic to $\SL(2,\R)$. Then $\mu$ is also ergodic with respect
to the one-parameter unipotent subgroup $U$ of $H$ corresponding to
the upper unipotent subgroup in $\SL(2,\R)$.
\end{corollary}

In fact, an invariant function $f \in L^2(X,\mu)$ that is invariant
under $U$ must be invariant under $\SL(2,\R)$ by Proposition
\ref{Mautner}. Since the latter group is assumed to be ergodic, the
function must be constant as required. (We leave it to the reader to
check the continuity requirement.)

\begin{proof}[Proof of Proposition \ref{Mautner}]
Following Margulis \cite{Margulis-Kyoto} we define the auxiliary
function $p (g) = (\phi (g) v, v)$. Notice first that the function
$p(\cdot)$ characterizes invariance in the sense that $p (g) = (v,
v)$ implies $\phi (g) v = v$. By continuity of the representation
$p(\cdot)$ is also continuous. Moreover, by our assumption on $v$
the map $p(\cdot)$ is bi-$U$-invariant since
\begin{equation*}
p(ugu')= (\phi (u) \phi (g) \phi (u')v,v)=(\phi (g) v, \phi (u ^{- 1}) v) = p (g).
\end{equation*}
Let $\epsilon, r,s \in \R$ and calculate
\begin{equation*}
\begin{pmatrix}
1 & r \\ 0 & 1
\end{pmatrix}
\begin{pmatrix}
1 & 0 \\ \epsilon & 1
\end{pmatrix}
\begin{pmatrix}
1 & s \\ 0 & 1
\end{pmatrix}=
\begin{pmatrix}
1+r \epsilon & r+s+rs \epsilon \\ \epsilon & 1+s \epsilon
\end{pmatrix}.
\end{equation*}
Now fix some $t \in \R$, let $\epsilon$ be close to zero but nonzero, choose $r=\frac{e^t-1}{ \epsilon}$ and $s=\frac{-r}{1+r \epsilon}$.  Then the above matrix simplifies to
\begin{equation*}
\begin{pmatrix} e^t & 0 \\ \epsilon & e ^{-t} \end{pmatrix}
\end{equation*}
In particular, this shows that
\begin{equation*}
p\left(\begin{pmatrix}
1 & 0 \\ \epsilon & 1
\end{pmatrix}\right)=p\left(\begin{pmatrix} e^t & 0 \\ \epsilon & e ^{-t} \end{pmatrix}
\right)\end{equation*}
is both close to $p(e)$ and to
\begin{equation*}
p\left(\begin{pmatrix} e^t & 0 \\ 0& e ^{-t} \end{pmatrix}\right).
\end{equation*}
Therefore, the latter equals $(v,v)$ which implies that $v$ is invariant under $\begin{pmatrix} e^t & 0 \\ 0 & e ^{-t} \end{pmatrix}$ as mentioned before.

The above implies now that $p(\cdot)$ is bi-invariant under the diagonal subgroup. Using this and the above argument once more, it follows that $v$ is also invariant under $\begin{pmatrix} 1 & 0 \\ s & 1 \end{pmatrix}$ for all $s \in \R$.
\end{proof}

\section{The proof of Theorem \ref{Main theorem}}

In this section we prove Theorem \ref{Main theorem} using the
prerequisites discussed in the last section. Let us mention again
that the general outline of the proof is very similar to the
strategy M.~Ratner \cite{Ratner-acta-measure} used to prove her
theorems.

From now on let $\mu$ be an $H$-invariant and ergodic probability
measure on $X = \Gamma \backslash G$.

\subsection{ The goal and the first steps}

It is easy to check that
\begin{equation*}
\operatorname{Stab}(\mu) =\{ g \in G:\mbox{ right multiplication with $ g$ on $X$ preserves }\mu\}
\end{equation*}
is a closed subgroup of $G$.  Let $L=\operatorname{Stab} (\mu) ^\circ$ be the connected component.  Then as discussed any element of $\operatorname{Stab} (\mu)$ sufficiently close to $e$ belongs to $L$.  Also since $\SL (2, \R)$ is connected we have $H <L$.

We will show that $\mu$ is concentrated on {\bf a single orbit} of
$L$, i.e. that there is some $L$-orbit $L.x_0$ of measure one $\mu
(L.x_0) = 1$. Then by $L$-invariance of $\mu$ and uniqueness of Haar
measure, $\mu$ would have to be the $L$-invariant volume form on
this orbit $L.x_0$. However, since $\mu$ is assumed to be a
probability measure this also implies that the orbit $L.x_0$ is
closed as seen in the next lemma.

\begin{lemma}\label{lemma to come}
If $\mu$ is concentrated on a single $L$-orbit $L.x_0$ and is $L$-invariant, then $L.x_0$ is closed and
$\mu$ is supported on $L.x_0$.
\end{lemma}

In the course of the proof we will recall a few facts about Haar
measures and also prove that a Lie group which admits a lattice is
unimodular, i.e.\ satisfies that the Haar measure is left and right
invariant. For a more comprehensive treatment of the relationship
between lattices and Haar measures see \cite{Raghunathan}.

\begin{proof}
Suppose $x _ i \in \ell _ i.x _ 0 \in L.x _ 0$ converges to $y$.  We have to show that $y \in L.x _ 0$.  Now either $x _ i \in \overline{B_1^L}. y$ for some $i$ --- i.e. the convergence is along $L$ and the lemma follows --- or $x _ i \not \in \overline{B_1^L}.y  $ for all $i$.  In the latter case we may choose a subsequence so that $x _ i \not \in \overline{B_1^L}.x _ j$ for $i \ne j$.  Now let $r < 1$ be an injectivity radius of $X$ at $y$.  Then $x _ i \in B _{r/2}.y$ for large enough $i$, say for $i \geq i _ 0$.  It follows that the sets $B_{r/2}^L.x_i \subset X$ are disjoint for $i \geq i _ 0$.  Since $x _ i = \ell _ i.x_0$ it follows that these sets are all of the form $B_{r/2}^L (\ell_i).x_0$.  We claim that the existence of the finite volume orbit implies that $L$ is unimodular.  If this is so, then we see that the sets $B_{r/2}^L.x_i$ all have the same measure, which contradicts the finite volume assumption.

It remains to show that $L$ is unimodular if the $L$-orbit of $x _ 0 = \Gamma g _ 0$ has finite $L $-invariant volume, or equivalently if $\Gamma _ L = g _ 0 ^{- 1} \Gamma g _ 0 \cap L < L$ is a lattice. So suppose $\mu$ is an $L$-invariant probability measure on $\Gamma _ L \backslash L$, where $L$ is acting on the right.  Then the following gives the relationship between $\mu$ and a right Haar measure $m_L$ on $L$.  Let $f$ be a compactly supported continuous function on $L$, then
\begin{equation*}
\int f \operatorname{d} \! m_L = \int \sum_{\gamma \in \Gamma _ L}f(\gamma \ell) \operatorname{d} \!  \mu.
\end{equation*}
Here note first that the sum $F (\ell)$ inside the integral on the right is finite for every $\ell$, satisfies $F(\gamma \ell) = F(\ell)$, and so defines a function on $\Gamma _ L \backslash L$ which is also compactly supported and continuous.  Therefore, the right integral is well-defined.  Using invariance of $\mu$ and the uniqueness of the right Haar measure on $L$, the equation follows (after possibly rescaling $m_L$).

The above formula immediately implies that $m_L$ is left-invariant under $\Gamma _ L$.  We use Poincar\'e recurrence to extend this to all of $L$.  Let $K \subset L$ be a compact subset of positive measure.  Let $\ell \in L$ be arbitrary and consider the map $T: \Gamma _ L \backslash L \rightarrow \Gamma _ L \backslash L$ defined by left multiplication by $\ell ^{- 1}$.  By assumption this preserves the probability measure $\mu$, therefore there exists some fixed $k \in K$, infinitely many $n _ i$, and $\gamma _{i} \in \Gamma _ L$ such that $\gamma _{ i} k \ell ^{- n _ i} \in K$.  In other words there exist infinitely many $k _ i \in K$ with $\ell ^{n _ i} = k _ i \nope ^{- 1} \gamma _ i k$.  The measure obtained by left multiplication by $\ell$ from a right Haar measure $m_L$ it again a right Haar measure and so must be a multiple $c (\ell)m_L$.  The constant $c (\ell) \in \R ^ \times$ defines a character, i.e.\ a continuous homomorphism $c: L \rightarrow \R ^ \times$.  It follows that $c (\ell)^{n _ i} = c(k _ i)c(\gamma _ i) c (k) = c (k _ i) c (k)$ remains bounded as $n_i \rightarrow \infty$.  Therefore, $c (\ell)=1$.  Since $\ell$ was arbitrary this proves the claim and the lemma.
\end{proof}

The main argument will be to show that if $\mu$ is not concentrated
on a single orbit of $L$, then there are other elements of
$\operatorname{Stab} (\mu)$ close to $e$.  This shows that we should
have started with a bigger subgroup $L'$. If we repeat the argument
with this bigger $L'$, we will either achieve our goal or make $L'$
even bigger.  We start by giving a local condition for a measure
$\mu$ to be concentrated on a single orbit.

\begin{lemma}  \label{Local condition}
Suppose $x _ 0 \in X$ has the property that $\mu (B_\delta^L.x_0)>0$ for some $\delta>0$, then $\mu$ is concentrated on $L.x_0$. So either the conclusion of Theorem \ref{Main theorem} holds for $L$ and $x_0$, or for every $x_0$ we have $\mu(B_\delta ^L.x_0)=0$.
\end{lemma}

\begin{proof}
This follows immediately from the definition of ergodicity of $\mu$ and the fact that $L.x_0$ is an $H$-invariant measurable set.
\end{proof}

We will be achieving the assumption to the last lemma by studying
large sets $X' \subset X$ of points with good properties.  Let $x _
0 \in X '$ be such that all balls around $x _ 0$ have positive
measure.  Suppose $X'$ has the property that points $x '$ close to
$x _ 0$ that also belong to $X '$ give `rise to additional
invariance' of $\mu$ unless $x$ and $x '$ are locally on the same
$L$-orbit (i.e.\ $x ' = \ell.x$ for some $\ell \in L$ close to $e$).
Then either $L$ can be made bigger or $B_\delta ^L.x\cap X'\subset
B_\delta^L.x_0$ for some $\delta>0$ and therefore the latter has
positive measure. However, to carry that argument through requires a
lot more work.  We start by a less ambitious statement where two
close by points in a special position from each other give `rise to
invariance' of $\mu$. Recall that $U=\left\{\begin{pmatrix} 1 &
* \\ 0 & 1\end{pmatrix}\right\}$.

\begin{proposition}  \label{Centralizer proposition}
There is a set $X ' \subset X$ of $\mu$-measure one such that if $x, x ' \in X '$ and $x ' = c. x$ with
\begin{equation*}
c \in C_G(U) =\{ g \in G: g u = u g\text{\ for all }u \in U\},
\end{equation*}
then $c$ preserves $\mu$.
\end{proposition}

The set $X'$ in the above proposition we define to be the set of $\mu$-generic points (for the one paramenter subgroup defined by $U$). A point $x \in X$ is {\em $\mu$-generic} if
\begin{equation*}
\frac{1}{T}\int_0^Tf(u_t.x)\operatorname{d}\!t \rightarrow \int f\operatorname{d}\!\mu\mbox{ for }T \rightarrow \infty
\end{equation*}
for all compactly supported, continuous functions $f:X \rightarrow \R$. Recall that by the Mautner phenomenon $\mu$ is $U$-ergodic. Now the ergodic theorem implies that the set $X '$ of all $\mu$-generic points has measure one.  (Here one first applies the ergodic theorem for a countable dense set of compactly supported, continuous functions and then extends the statement to all such functions by approximation.)

\begin{proof}
For $c \in C _ G (U)$ and a compactly supported, continuous function $f: X \rightarrow \R$ define the function $f_c(x)=f(c.x)$ of the same type.  Now assume $x,x'=c.x \in X'$ are $ \mu$-generic.  Since $u _ t c =cu _ t$ we have $f(u _ t.x') = f (cu _ t.x) = f_c (u _ t.x)$ and so the limits
\begin{align*}
\frac{1}{T}\int_0^Tf(u_t.x')\operatorname{d}\!t & \rightarrow \int f\operatorname{d}\!\mu\mbox{ and}\\
\frac{1}{T}\int_0^Tf_c(u_t.x)\operatorname{d}\!t & \rightarrow \int f_c\operatorname{d}\!\mu
\end{align*}
are equal.  However, the last integral equals $\int f_c\operatorname{d}\!\mu=\int f\operatorname{d}\!c_*\mu$
were $c_*\mu$ is the push forward of $\mu$ under $c$. Since $f$ was any compactly supported, continuous function, $\mu = c_*\mu$ as claimed.
\end{proof}

\subsection{Outline of the H-principle}
In Proposition \ref{Centralizer proposition} we derived invariance of $\mu$ but only if we have two points $x, x ' \in X '$ that are in a very special relationship to each other.  On the other hand if $\mu$ is not supported on the single $L$-orbit, then we know that we can find many $y, y ' \in X'$ that are close together but are not on the same $L$-leaf locally by Lemma \ref{Local condition}.  Without too much work we will see that we can assume
\begin{equation*}
y ' = \exp (v).  y \mbox{ with }v \in \mathfrak l'
\end{equation*}
where $ \mathfrak l'$ is an $\SL (2, \R)$-invariant complement in
$\mathfrak g$ of the Lie algebra $\mathfrak l$ of $L$, see Lemma
\ref{pigeon hole}. What we are going to describe is a version of the
so-called H-principle as introduced by Ratner
\cite{Ratner-factor,Ratner-joinings} and generalized by Witte
\cite{Witte-rigidity}, see also \cite{Witte}.

By applying the same unipotent matrix $u \in U$ to $y$ and $y '$ we get
\begin{equation*}
u.  y ' = (u \exp (v) u ^{- 1}).(u.y) = \exp (\operatorname{Ad} _ u (v)).(u.y).
\end{equation*}
In other words, the divergence of the orbits through $y$ and $y '$ can be described by conjugation in $G$ --- or even by the adjoint representation on $\mathfrak g$.  Since $H$ is assume to be isomorphic to $\SL (2, \R)$ we will be able to use the theory on representations as in Section \ref{Representations}.
In particular, recall that the fastest divergence is happening along a direction which is stabilized by $U$.
Since all points on the orbit of a $\mu$-generic point are also $\mu$-generic, one could hope to flow along $U$ until the two points $x=u.y, x '=u.x'$ differ significantly but not yet to much.  Then $y'=u.x'= h.(u.x)=h.y$ with $h$ almost in $C_G(U)$.  To fix the `almost' in this statement we will consider points that are even closer to each other, flow along $U$ for a longer time, and get a sequence of pairs of $\mu$-generic points that differ more and more by some element of $C_G(U)$.  In the limit we hope to get to points that differ precisely by some element of $C_G(U)$ which is not in $L$.

The main problem is that limits of $\mu$-generic points need not be $\mu$-generic (even for actions of unipotent groups).  Therefore, we need to introduce quite early in the argument a compact subset $K \subset X '$ of almost full measure that consists entirely of $\mu$-generic points.  When constructing $u.x',u.x$ we will make sure that they belong to $K$ --- this way we will be able to go to the limit and get $\mu$-generic points that differ by some element of $C_G(U)$.

We are now ready to proceed more rigorously.

\subsection{Formal preparations, the sets $K$, $X _ 1$, and $X _ 2$}

Let $X '$ be the sets of $\mu$-generic points as above, and let $K \subset X '$ be compact with $\mu (K) > 0.9$.  By the ergodic theorem
\begin{equation*}
\frac{1}{T} \int_ 0 ^ T 1 _ K (u _ t.y) \operatorname{d}\!t \rightarrow \mu (K)
\end{equation*}
for $\mu$-a.e.\ $y \in X$.  In particular, we must have for a.e.\ $y \in X$
\begin{equation*}
\frac{1}{T} \int_ 0 ^ T 1 _ K (u _ t.y) \operatorname{d}\!t > 0.8\mbox{ for large enough } T.
\end{equation*}
Here $T$ may depend on $y$ but by choosing $T _ 0$ large enough we may assume that the set
\begin{equation*}
X_1= \left \{ y \in X:\frac{1}{T} \int_ 0 ^ T 1 _ K (u _ t.y) \operatorname{d}\!t > 0.8 \text{{ for all }} T\geq T_0 \right\}
\end{equation*}
has measure $\mu (X _ 1) >0.99$.  By definition points in $X _ 1$ visit $K$ often enough so that we will be able to find for any $y, y ' \in X_1$ many common values of $t$ with $u _ t.y, u _ t.y'\in K$.

The last preparation we need will allow us to find $y, y ' \in X _ 1$ that differ by some $\exp (v)$ with $v \in \mathfrak l '$. For this we define
\begin{equation*}
X _ 2 = \left \{ z \in X:\frac{1}{m_L(B_1^L)} \int_ {B_1^L} 1 _ {X _ 1} (\ell.  z) \operatorname{d} \!  m_L (\ell) > 0.9\right \}
\end{equation*}
where $ m _ L$ is a Haar measure on $L$. (Any other smooth measure would do here as well.)

\begin{lemma} $\mu (X _ 2) > 0 .9$.  \end{lemma}

\begin{proof}
Define $Y = X \times B_1^L$ and consider the product measure $\nu = \frac{1}{m_L(B_1^L)} \mu \times m_L$ on $Y$.  The function $f(z,\ell)=1_{X_1}(\ell.z)$ integrated over $z$ gives independently of $\ell$ allways $\mu (X _ 1)$ since $\mu$ is $L$-invariant.  By integrating over $\ell$ first we get therefore the function
\begin{equation*}
g(z)=\frac{1}{m_L(B_1^L)} \int_ {B_1^L} 1 _ {X _ 1} (\ell.  z) \operatorname{d} \!  m_L (\ell)
\end{equation*}
whose integral satisfies $\int g \operatorname{d} \!  \mu = \mu (X _ 1)$.  Since $g (z) \in [0, 1]$ for all $z$
\begin{equation*}
0 .99 < \int_{X _ 2} g \operatorname{d} \!  \mu+\int_{X\setminus X_2} g \operatorname{d} \!  \mu \leq \mu (X _ 2) + 0 .9 \mu (X \setminus X_2) = 0 .1 \mu (X _ 2) + 0 .9
\end{equation*}
which implies the lemma.
\end{proof}

Let as before $\mathfrak l \subset \mathfrak g$ be the Lie algebra of $L < G$ and let $\mathfrak l ' \subset  \mathfrak g$ be an $\SL (2, \R)$-invariant complement of $\mathfrak l $ in $\mathfrak g$.  Then the map $\phi: \mathfrak l ' \times \mathfrak l  \rightarrow G$ defined by $\phi ( v,w) = \exp (v) \exp (w)$ is $C^ \infty$ and its derivative at $(0,0)$ is the embedding of $\mathfrak l ' \times \mathfrak l $ into $\mathfrak g$.  Therefore, $\phi$ is locally invertible so that every $g \in G$ close to $e$ is a unique product $g = \exp (v)\ell$ for some $\ell \in L$ close to $e$ and some small $v \in \mathfrak l '$.  We define $\pi _ L (g) = \ell$.  For simplicity of notation we assume that this map is defined on an open set containing $\overline{B_1 ^ L}$ (if necessary we rescale the metric).

\begin{lemma}\label{pigeon hole}
For any $\epsilon > 0$ there exists $\delta > 0$ such that for $g \in B_\delta ^ G$, and $z, z ' = g.z \in X _ 2 $ there are $\ell_2 \in B_1^L$ and $\ell_1 \in B_{ \epsilon} ^ L(\ell_2)$ with $\ell_1 .z, \ell_2.z'\in X _ 1$ and $\ell_2g\ell_1^{-1}= \exp (v)$ for some $v \in B_\epsilon ^{ \mathfrak l '} (0)$.
\end{lemma}

\begin{proof}
Let $g \in B_\delta ^ G$ and consider the $C ^ \infty$-function $\psi (\ell) = \pi _ L (\ell g)$.  If $g = e$ then $\psi$ is the identity, therefore if $\delta$ is small enough the derivative of $\psi$ is close to the identity and its Jacobian is close to one. In particular, we can ensure that $m_L(\psi (E ' )) > 0.9 m_L(E ')$ for any measurable subset $E ' \subset B_1^L$.  Moreover, for $\delta$ small enough we have $\psi (\ell) \in B_{\epsilon} ^ L(\ell)$ for any $\ell \in B_1 ^ L$ and so  $\psi (B_1^L) \subset B_{ 1 + \epsilon} ^ L$.

Now define the sets $E = \{\ell \in B _ 1 ^ L: \ell.z \in X _ 1 \}$ and $E ' = \{\ell \in B _ 1 ^ L: \ell.z ' \in X _ 1 \}$ which satisfy $m_L(E)>0.9m_L(B_1^L)$ and $m_L(E ')>0.9m_L(B_1^L)$ by definition of $X _ 2$. We may assume that $\epsilon$ is small enough so that $m_L(B_{ 1 + \epsilon}^L)<1.1 m_L(B_1^L)$. Together with the above estimate we now have $m_L(\psi (E')) > 0.5 m _ L (B_{1 + \epsilon} ^L)$ and $m_L(E) > 0.5 m _ L (B_{1 + \epsilon} ^L)$.  Therefore, there exists some $\ell_2 \in E'$ with $\ell_1= \psi (\ell_2) \in E$.  By definition of $E,E'$ we have $\ell_1.y,\ell_2.y ' \in X _ 1$.  Finally, by definition of $\pi _ L$
\begin{equation*}
\ell_2 g \ell_1 ^ { - 1} =\ell_2 g \pi _ L (\ell_2 g) ^{- 1} = \exp (v)
\end{equation*}
for some $v \in \mathfrak l '$. Again for sufficiently small $\delta$ we will have $v \in B_\epsilon ^{ \mathfrak l '} (0)$.
\end{proof}

\subsection{H-principle for $\SL(2,\R)$}

Let $x_0 \in X_2 \cap \supp \mu|_{X_2}$ so that
\begin{equation*}
\mu\bigl((B_\delta^G.x_0) \cap X_2\bigr)>0 \text {  for all } \delta>0.
\end{equation*}
Now one of the following two statements must hold:
\begin{enumerate}
\item
there exists some $\delta>0$ such that $B_\delta^G.x_0 \cap X_2 \subset B_\delta^L.x_0$, or
\item
for all $\delta>0$ we have $B_\delta^G.x_0 \cap X_2 \not \subset B_\delta ^L.x_0$.
\end{enumerate}

We claim that actually only (1) above is possible if $L$ is really the connected component of $\operatorname{Stab}(\mu)$. Assuming this has been shown, then we have $\mu(B_\delta ^L.x_0)>0$ which was the assumption to Lemma \ref{Local condition}. Therefore, $\mu(L.x_0)=1$ and by Lemma \ref{lemma to come} $L.x_0 \subset X$ is closed --- Theorem~\ref{Main theorem} follows. So what we really have to show is that (2) implies that $\mu$ is invariant under a one parameter subgroup that does not belong to $L$.

\begin{lemma}  \label{Transverse existence}
Assuming (2) there are for every $\epsilon >0$ two points $y, y '
\in X _ 1$ with $d (y, y ') < \epsilon$ and $y ' = \exp (v). y$ for
some nonzero $v \in B_\epsilon^{\mathfrak l'}(0)$.
\end{lemma}

\begin{proof}
Let $z=x_0$. By (2) there exists a point $z' \in X_2$ with $d(z,z')<\delta$ and $z' \not \in B_{\delta}^L.z$. Let $g \in B_{\delta}^G$ be such that $z'=g.z$ and $g \not \in L$. Applying Lemma~\ref{pigeon hole} the statement follows since $g \not \in L$ and so $v \ne 0$ by our choice of $z'$.
\end{proof}

Using $y,y' \in X_1$ and $v \in \mathfrak l'$ for all $\epsilon>0$ as in the above lemma we will show that $\mu$ is invariant under a one-parameter subgroup that does not belong to $L$. For this it is enough to show the following:

{\bf Claim:} For any $\eta>0$ there exists a nonzero $w \in B_\eta^{\mathfrak l'}(0)$ such that $\mu$ is invariant under $\exp(w)$.

To see that this is the remaining assertion, notice that we then also have invariance of $\mu$ under the subgroup $\exp(\Z w)$. While this subgroup could still be discrete, when $\eta \rightarrow 0$ we find by compactness of the unit ball in $\mathfrak l'$ a limiting one parameter subgroup $\exp(\R w)$ that leaves $\mu$ invariant.

We start proving the claim. Let $\eta>0$ be fixed, and let
$\epsilon>0$, $y,y' \in X_1$, and $v \in B_\epsilon ^{\mathfrak
l'}(0)$ as above. We will think of $\epsilon$ as much smaller than
$\eta$ since we will below let $\epsilon$ shrink to zero while not
changing $\eta$. Let $\operatorname{Sym}^n(\R ^2)$ be an irreducible
representation as in Section \ref{Representations}, and let
$p=p(A,B) \in \operatorname{Sym}^n(\R ^2)$. Recall that
$u_t=\begin{pmatrix} 1 & t \\ 0 & 1\end{pmatrix}$ applied to
$p(A,B)$ gives $p(A,B+tA)$. We define
\[
T_p=\frac{\eta}{\max(|c_1|,\ldots,|c_n|^{1/n})}
\]
and set $T_p=\infty$ if the expression on the right is not defined. The significance of $T_p$ is that for $t=T_p$ at least one term in the sum $(c_0+c_1t+\cdots+c_nt^n)$ is of absolute value one while all others are less than that --- recall that this sum is the coefficient of $A$ in $p(A,B+tA)$.
To extend this definition to $\mathfrak l'$ which is not necessarily irreducible we split $\mathfrak l'$ into irreducible representations $\mathfrak l'=\bigoplus_{j=1}^k V_j$ and define for $v=(p_j)_{j=1,\ldots,k}$
\[
T_v=\min_j T_{p_j}.
\]

\begin{lemma}  \label{Next to last}
There exists constants $n>0$ and $C>0$ that only depend on
$\mathfrak l'$ such that for  $v \in B_\epsilon ^{\mathfrak l'}(0)$
and $t \in [0,T_v]$ we have
\begin{equation*}
\operatorname{Ad} _{u_t}(v)
= w+O(\epsilon ^{1/n})
\end{equation*}
where $w \in B_{C \eta}^{\mathfrak l'}(0)$ is fixed under the subgroup $U=u_{\R}$. Here we write $O(\epsilon ^{1/n})$ to indicate a vector in $\mathfrak l'$ of norm less than $C \epsilon^{1/n}$.
\end{lemma}

\begin{proof}
We first show the statement for irreducible representations $\operatorname{Sym} ^ n (\R ^ 2)$ inside $\mathfrak l'$.  There any multiple of $A ^ n$ is fixed under $U$.  Similar to the earlier discussion the coefficient $(c _ 0 + c _ 1 t + \cdots + c _ n t ^ n)$ is for $t \in [0, T _ p]$ bounded by $n \eta$.  For the other coefficients note first that these are sums of terms of the form $c _ i t ^ j$ for $j < i$.  Since $t < \eta | c _ i| ^{- 1 / i}$ we have that each such term is bounded by $|c _ i t ^ j| \leq | c _ i|| c _ i| ^{-j/i} \leq| c _ i| ^{1 / i} \ll \epsilon ^{1/n}$ where the implied constant only depends on the norm on $\mathfrak l '$ and the way we split $\mathfrak l '$ into irreducibles (and we assumed $\eta<1$).  This proves the lemma for irreducible representations.

The general case is now straightforward.  If $t < T _ v = \min_ j T _ {p _ j}$ then since $\operatorname{Ad} _{u _ t} (p _ j)$ is of the required form for $j = 1,\ldots, k$ the lemma follows by taking sums.
\end{proof}

If $v$ is already fixed by $U$ then $T _ v = \infty$ (and other way around) and the above statement is rather trivial since $w = v$.  Moreover, by definition of $X _ 1$ we have $\frac{1}{T} \int_ 0 ^ T 1_K(u _ t x _ i) \operatorname{d} \! t >0.8 $ for $i = 1, 2$.  From this it follows that there is some $t \in [0, T _ 0]$ with $u _ t.x _ 1, u _ t.x _ 2 \in K$.  Since $K \subset X '$ Lemma \ref{Centralizer proposition} proves (assuming $\epsilon < \eta$) the claim in that case and we may from now on assume that $v$ is not fixed under the action of $U$ and so $ T _ v < \infty$.

\begin{lemma}  \label{Last one}
There exists a constant $c > 0$ that only depends on $\mathfrak l '$ such that the decomposition $\operatorname{Ad} _{u _ t} (v) = w + O(\epsilon ^{1/n})$ as above satisfies $\| w\| >c \eta$ for $t \in E_v$ where $E_v \subset [0,T_v]$ has Lebesgue measure at least $0 .9 T _ v$.
\end{lemma}

\begin{proof}
We only have to look at the irreducible representations $V_j= \operatorname{Sym} ^ n (\R ^ 2)$ in $\mathfrak l '$ with $T _ v = T _{ p _ j}$.  The size of the corresponding component of $w$ is determined by the value of the polynomial $ c _ 0 + c _ 1 t +\cdots+ c _ n t ^ n$.  We change our variable by setting $s=\frac{t}{ T _ v}$ and get the polynomial $q (s) = c _ 0 +\frac{ c _ 1}{ T _ v} s +\cdots +\frac{ c _ n}{ T _ v ^ n} s ^ n$. By definition of $T _ v$ the polynomial $q(s)$ has at least one coefficient of absolute value one while the others are in absolute value less or equal than one. Therefore, we have reduced the problem to finding a constant $c$ such that for every such polynomial we have that $E_1=\{s \in [0,1]:|q(s)|\geq c\}$ has measure bigger than $0.9$.

This can be done in various ways --- we will use the following property of polynomials for the proof. Every polynomial of degree $n$ is determined by $n+1$ values (by the standard interpolation procedure). Moreover, we can give an upper bound on the coefficients in terms of these values unless the values of $s$ used for the data points are very close together. (The determinant of the Vandermonde determinant is the product of the differences of the values of $s$ used in the interpolation.) Suppose $s_0  \in [0,1]\setminus E_1$, then all points close to $s_0$ give unsuprisingly also small values. So we now look for $s_1 \in [0,1] \setminus (E_1 \cup [s_n-\frac{1}{20n},s_n+\frac{1}{20n}])$ --- if there is no such point then $E_1$ is as big as required. We repeat this search until we have found $n+1$ points $s_0,\ldots,s_n \in [0,1]\setminus E_1$ that are all separated from each other by at least $\frac{1}{20n}$. Again if we are not able to find these points, then $E_1$ is sufficiently big. However, as discussed the coefficients of $q$ are due to $|q(s_0)|,\ldots,|q(s_n)|<c$ bounded by a multiple of $c$ (that also involves some power of the degree $n$). If $c$ is small enough compared to that coefficient we get a contradiction to the assumption that $q$ has at least one coefficient of absolute value one. It follows that for that choice of $c$ we can find at most $n$ points in the above search and so $E_1$ has Lebesgue measure at least $1-n\frac{2}{20n}=0.9$.
\end{proof}

Recall that case (1) from the beginning of this section implies Theorem \ref{Main theorem} and that we are assuming case (2).  Moreover, recall that this implies for all $\epsilon > 0$ the existence of $y, y ' \in X _ 1$ with $y ' = \exp (v).  y$ for some nonzero $v \in B_\epsilon ^{\mathfrak l '} (0)$ by Lemma \ref{Transverse existence}.  By definition of $X _ 1$ the sets
\begin{align*}
E_T & = \{ t \in [0, T]: u _ t . y \in K \} \text{and} \\
E_T ' & = \{ t \in [0, T]: u _ t . y ' \in K \}
\end{align*}
have Lebesgue measure bigger than $0.8 T$ whenever $T \geq T _ 0$. From the definition it is easy to see that $T_v \geq T_0$ once $\epsilon$ and therefore $v$ are sufficiently small, so we can set to $T = T _ v$.  Moreover, let $E_v$ be as in Lemma \ref{Last one}.  Then the union of the complements of these three sets in $[0, T _ v]$ has Lebesgue measure less than $0 .5 T_v$. Therefore, there exists some $t \in E_{ T _ v} \cap E_{ T _ v} ' \cap E _ v$.  We set $x = u _ t. y$, $x ' = u _ t. y '$ which both belong to $K$ by definition of $E _{T _ v}$ and $E _{T _ v} '$.  Moreover, $x ' = \exp (w + O(\epsilon ^{1/n})).  X$ where $ w \in \mathfrak l'$ is stabilized by $U$ and satisfies $c \eta \leq \|w\| \leq C \eta$ by Lemma~\ref{Next to last}-- \ref{Last one}.  We let $\epsilon \rightarrow 0$ and choose converging subsequences for $x, x '$, and $w$.  This shows the existence of $x, x ' = \exp (w). x \in K$ and $w \in \mathfrak l '$ with $c \eta \leq\| w\| \leq C \eta$ which is stabilized by $U$. This is in effect our earlier claim which as we have shown implies that $\mu$ is invariant under a oneparameter subgroup not belong into $L$. This concludes the proof of Theorem \ref{Main theorem}.



\def\cprime{$'$}

\end{document}